\documentclass[12pt]{article}

\usepackage[utf8]{inputenc}

\usepackage{tikz}
\usepackage{amsmath} 
\usepackage{bm}
\usepackage{multirow}
\usepackage{booktabs}

\usepackage{hyperref}

\graphicspath{ {./figs/}}

\newcommand{\grad}{\mathop{\rm grad}\nolimits}
\renewcommand{\div}{\mathop{\rm div}\nolimits}

\makeatletter
\renewcommand{\@biblabel}[1]{#1.\hfill}
\makeatother

\begin{document}
\title{Meshfree multiscale method with partially explicit time discretization}

\author{Nikiforov Djulustan\\ North-Eastern Federal University, Yakutsk, Russia}

\date{}
\maketitle
\begin{abstract} 
In this paper, a multiscale approach with partially explicit time discretization is proposed. The idea is to use a partially explicit time scheme, considering a filtration problem in a fractured medium, where the implicit scheme is used for nodes whose subdomains contain fractures, and the explicit scheme is used for all others. In this way, it is possible to use a time step that is independent of the diffusion coefficient for fractures. Numerical results demonstrating high accuracy of calculations are presented.
\end{abstract}

\section{Introduction}
Multiscale methods are relevant in various fields where phenomena occur at multiple scales or levels of detail. In numerical simulations and modeling, multiscale methods are used to address problems that involve multiple spatial or temporal scales. By considering different scales simultaneously, these methods can provide a more comprehensive understanding of the system's behavior and capture important features that may be missed by traditional numerical methods. Overall, multiscale methods play an important role in bridging the gap between different scales of analysis, allowing for a more comprehensive understanding of complex systems \cite{Efendiev1, Efendiev2}.

The contrast of parameters refers to significant variations or disparities in the properties or characteristics of the system being studied across different scales. This can include variations in material properties, boundary conditions, or other relevant parameters. The contrast of parameters strongly influences the solution of a multiscale problem. Multiscale methods are specifically designed to solve these problems by constructing basis functions that can efficiently handle parameter variations across different scales \cite{Nikiforov1, Ammosov2}.

Nevertheless, the high contrast of parameters can affect the choice of time integration schemes in time-dependent problems. Explicit time schemes rely on small time steps to ensure stability and accuracy. However, in the presence of high parameter contrast, the time step size may need to be further reduced to capture the fast dynamics associated with the high-contrast regions. This can significantly increase the computational cost and make the explicit time scheme impractical or inefficient for such problems. 

In such cases, alternative time integration schemes, such as implicit or partially explicit methods, may be more suitable to handle the high contrast and ensure stability and efficiency \cite{VPN1}. For example, \cite{Ammosov1} proposes a new approach to solve the problem of poroelasticity in fractured media based on hybrid explicit-implicit learning (HEI). To facilitate the calculations, the authors use the fixed strain and fixed stress splitting schemes. The main idea of the proposed method lies in partial learning. The proposed splitting algorithm in \cite{Efendiev3} treats dominant multiscale modes in the implicit fashion, while the rest in the explicit fashion. The contrast-independent stability of these algorithms requires a special multiscale space design, which is the main purpose of the paper. It is shown that with an appropriate choice of multiscale spaces, an unconditional contrast stability can be achieved.

In this paper, we propose a multiscale approach with partially explicit time discretization. This approach is based on a meshfree generalized multiscale finite element method \cite{Nikiforov3, Nikiforov4}. The idea of the splitting scheme is to split the vector of coarse mesh solutions for the stiffness matrix into two parts: with fast processes (with fractures) and slow processes (without fractures). Next, use an implicit time scheme for nodes with fast processes, and an explicit time scheme for all others. In this way, it is possible to use a time step that is independent of the diffusion coefficient for fractures.

\section{Problem formulation}

In this paper, we consider the problem of single-phase filtration in a fractured porous medium. The mathematical model is described by a coupled system of differential equations for the pressure $p_m$ in the porous matrix $\Omega$ and the pressure $p_f$ in the fractures $\gamma$   
\begin{equation}\label{eq:equation1}\begin{split}
c_m \frac{\partial p_m}{\partial t} - \div \left( \frac{k_m}{\mu} \grad p_m \right) + q(p_m-p_f)= 0, \quad x \in \Omega, \\
c_f \frac{\partial p_f}{\partial t} - \div \left( \frac{k_f}{\mu} \grad p_f \right) + q(p_f-p_m)= 0, \quad x \in \gamma, 
\end{split}\end{equation} 
where $c_i$ is the media compressibility, $k_i$ is the media permeability, $\mu$ is the fluid viscosity and $q$ is the mass transfer coefficient ($i = m, f$). 

System of equation \eqref{eq:equation1} is supplemented by the initial and boundary conditions
\begin{equation}\begin{split}
p_i = p_0, \quad &t = 0,    
\\
p_m = g, \quad &x \in \xi,
\\
p_f = g, \quad &x \in \xi,
\\
-\frac{k_m}{\mu}(\grad p_m, n) = 0, \quad &x \in \partial \Omega,
\\
-\frac{k_f}{\mu}(\grad p_f, n) = 0, \quad &x \in \partial \gamma.
\nonumber
\end{split}\end{equation} 

\begin{figure}[h!]
\begin{center}
\center{\includegraphics[width=\linewidth]{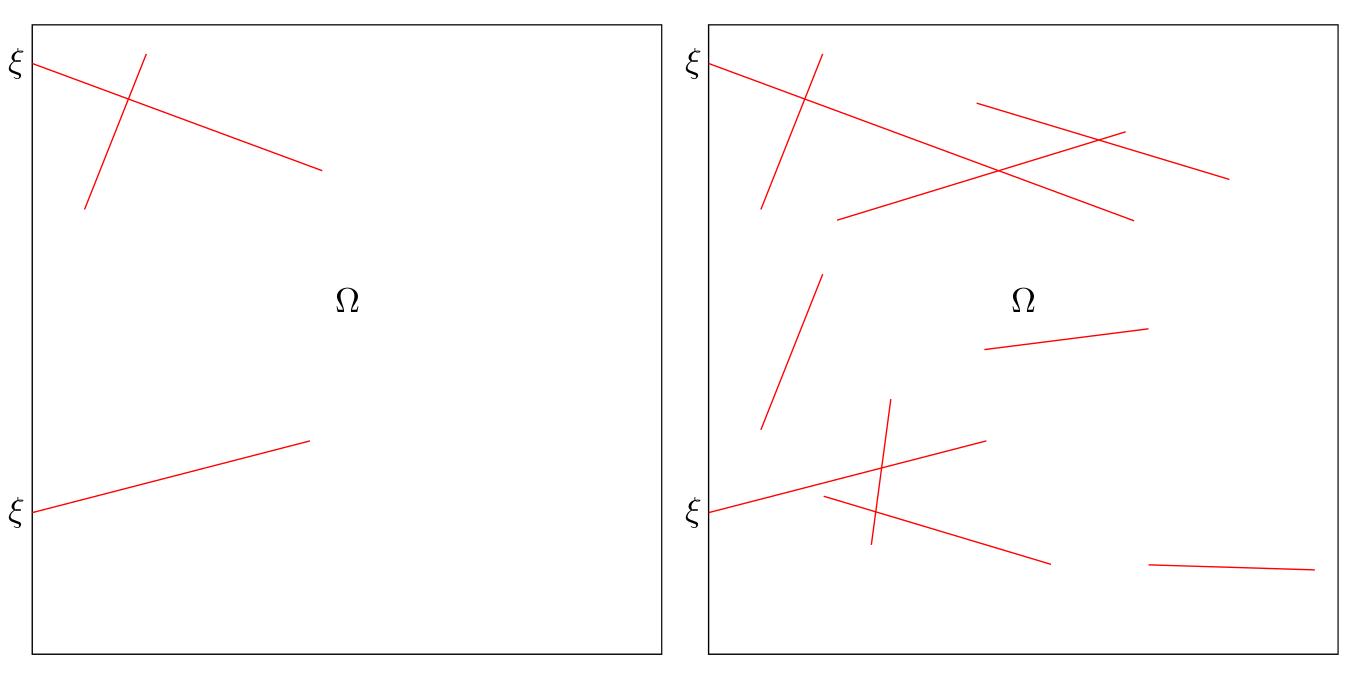}}
\end{center}
\caption{Computational domains $\Omega$ with the fractures $\gamma$ (red lines). Left: Test 1. Right: Test 2.}
\label{ris:domain}
\end{figure}

For our numerical tests, we take two different computational domains as shown in the Figure \ref{ris:domain} to show the robustness of our proposed approach. Here, in the first geometry there are few fractures (Test 1), and in the second there is a more complicated permeability field with a large number of fractures (Test 2).  

\section{Fine-grid discretization}

Let us define an uniform, for simplicity, time grid $\omega_t = \{ t^n = n \tau, \ n = 0, 1, \ldots, N_t-1, \ \tau N_t = t_{max} \}$, where $\tau$ - time step, $N_t$ - number of time steps and $t_{max}$ - solution time. Here we use the implicit time discretization of the system \eqref{eq:equation1}. For spatial approximation on a fine grid, we use the continuous Galerkin finite element method with standard linear basis functions (FEM). For the fractures accounting we use the discrete fracture model (DFM) \cite{dfm1, dfm2}. We use the superposition principle and assume that $p_1 = p_2 = p$. Then, summing the equations of the system \eqref{eq:equation1}, we obtain the following variational problem: find $p^{n+1} \in V$ such that 
\begin{equation}
m(\frac{p^{n+1}-p^{n}}{\tau}, v) + a(p^{n+1}, v) = 0, \quad n = 0, 1, \ldots, N_t-1,
\nonumber
\end{equation} 
where
\begin{equation}\label{eq:app2}\begin{split}
m = m_1 + m_2, \quad a = a_1 + a_2, \\
m_1(\frac{p^{n+1}-p^{n}}{\tau}, v) = \int\limits_{\Omega}c_m\frac{p^{n+1}-p^{n}}{\tau}vdx, \\
m_2(\frac{p^{n+1}-p^{n}}{\tau}, v) = \alpha\int\limits_{\gamma}c_f\frac{p^{n+1}-p^{n}}{\tau}vds, \\
a_1(p^{n+1}, v) = \int\limits_{\Omega} \left( \frac{k_m}{\mu} \grad p^{n+1}, \grad v \right)dx, \\
a_1(p^{n+1}, v) = \alpha\int\limits_{\gamma} \left( \frac{k_f}{\mu} \grad p^{n+1}, \grad v \right)ds, 
\end{split}\end{equation}
where $\forall v \in \hat{V}$. Here $p^{n} \approx p(t^n)$ and $\alpha$ is the fracture aperture.

We build an unstructured mesh $\mathcal{T}_h$ and construct discrete function space $V_h \subset V$. Then we use $p_h = (p_{h, 1}, p_{h, 2}, \ldots, p_{h, Nf})^T$ to denote the vector of the required unknowns, where $N_f$ -- the number of fine mesh vertices. Next we write the following matrix form for the fully discrete system 
\begin{equation}\label{eq:app3}
M \frac{p_h^{n+1}-p_h^{n}}{\tau} + A p_h^{n+1} = 0,
\end{equation}
where
\begin{equation}\label{eq:app4}\begin{split}
M = \{m_{ij}\}, \quad m_{ij} = m_1(\phi_j, \phi_i) + m_2(\hat{\phi}_j, \hat{\phi}_i), \\ 
A = \{a_{ij}\}, \quad a_{ij} = a_1(\phi_j, \phi_i) + a_2(\hat{\phi}_j, \hat{\phi}_i), 
\end{split}\end{equation}
where $\phi_i$ -- two-dimensional and $\hat{\phi}_i$ -- one-dimensional fine-scale piecewise linear basis functions for the porous medium and fractures.

\section{Meshfree GMsFEM}

In this section, we give a brief description of the Meshfree GMsFEM algorithm for our problem. In offline stage we have the following steps
\begin{enumerate}
\item Generation of meshfree coarse grid depending on fractures.
\item Construction of a multiscale offline space.
\item Assembling a projection matrix into a multiscale offline space.
\end{enumerate}
In the online stage, the fine-grid system is converted to a coarse grid system using the resulting projection matrix. As a result, the problem \eqref{eq:app2} is solved on a coarse grid.

\subsection{Meshfree coarse scale}
\begin{figure}[h!]
\begin{center}
\includegraphics[width=0.5\linewidth]{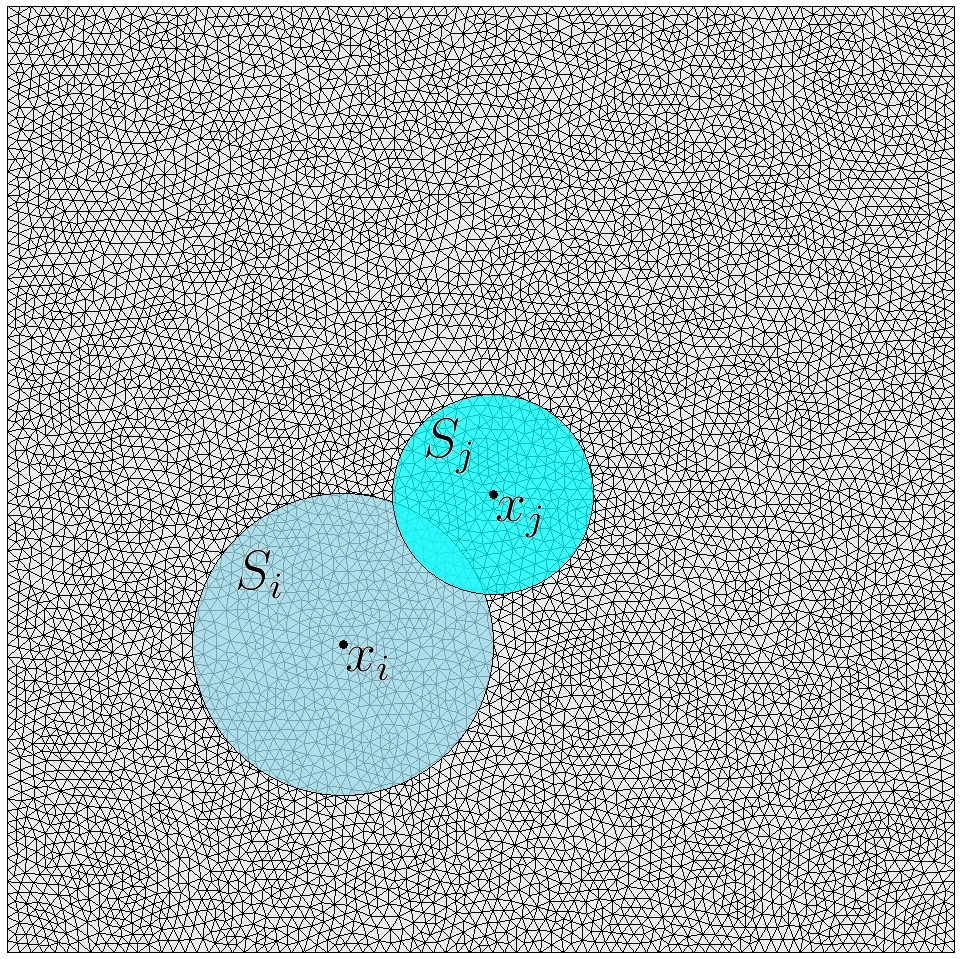}
\end{center}
\caption{Example of some coarse-scale elements $S_i$ and $S_j$.}
\label{ris:coarse_example}
\end{figure}

In Meshfree GMsFEM we use a point cloud instead of a structured triangular coarse grid. Let $\mathcal{S}_H$ be a partition of the computational domain $\Omega$ to the point cloud so that $\Omega \subset \bigcup_{i=1}^{N} S_i$ and suppose that each coarse element $S_i$ is partitioned into a connected union of elements of a fine grid (Fig. \ref{ris:coarse_example}). Let ${\{ x_i \}}^{N}_{i=1}$ is the coarse-scale nodes, where $N$ denotes the number of coarse nodes. Here, the coarse elements $S_i$ are the supports of basis functions
\begin{equation}
S_i = \{ y \in \mathcal{R}^N : \ \parallel y - x_i \parallel \le r_i \},
\nonumber
\end{equation} 
where $r_i$ is radius of coarse element $S_i$.

We denote the basis functions by $\psi_{i, k}$, which is supported in $S_i$, and the index $k$ represents the numbering of these basis functions. In turn, the solution will be sought as
\begin{equation}\label{eq:ms_sol}
p_{\text{c}}(x) = \sum_{i, k} c_{i, k} \psi_{i, k}(x).
\end{equation} 
Here the function $p_{\text{c}}(x)$ is general, referring to the pressure. This is due to the fact that the problem for the pressure is solved sequentially and for problem \eqref{eq:app2} the multiscale method will be almost the same. Where the multiscale method differs, we will point it out.

The main difficulties of the meshfree method are the choice of the location of points $\{x_i\}_{i=1}^{N}$, as well as the sizes of $r_i$ for $\{S_i\}_{i=1}^{N}$. The version of CVT (centroidal Voronoi tessellations) proposed in \cite{cvt1} is used in this work. An important point here is the choice of the distribution density function $\rho(x)$. Here we calculate the distribution density function $\rho(x)$ of a random variable taking into account the fractures using the following problem
\begin{equation}\label{eq:pc_problem}
\begin{split}
- \div \left( \beta \grad \rho \right) + \rho = f, \quad &x \in \Omega\\
- \beta \frac{\partial \rho}{\partial {n} } = 0, \quad &x \in \partial \Omega,
\end{split}
\end{equation}
where $f$ is a piecewise linear function depending on the location of the fractures and $\beta$ is a smoothing parameter, which in our case $\beta = 5$ was chosen by numerical experiments. Here the right hand side in discrete form is defined as 
\begin{equation}
f_h(x_i) =  
\begin{cases} 
10^5, \ &x_i \in \gamma, \\ 
1, \ &\text{otherwise},
\end{cases}
\nonumber
\end{equation}
where $i = 1, 2, \ldots, N_f$ and $N_f$ is a number of nodes of the fine grid. Next, we use the algorithm proposed in \cite{cvt2} for finding the radii $r_i$ for $S_i$. 

\subsection{Local basis functions}
To construct the offline space $V^{S}_{\text{off}}$, we solve the following local eigenvalue problem in each $S_i$
\begin{equation}\label{eq:spectral_problem}
A \Psi^{\text{off}} = \lambda^{\text{off}} B \Psi^{\text{off}},
\nonumber
\end{equation}
where $A$ and $B$ denote similar fine-scale matrices as
\begin{equation}
\begin{split}
A = \{a_{ij}\}, \quad a_{ij} = \int_{S} \left(\lambda_1 \grad \phi_j, \grad \phi_i \right) dx + \alpha \int_{\gamma^S} \left( \lambda_2 \grad \hat{\phi}_j, \grad \hat{\phi}_i \right) ds,\\ 
B = \{b_{ij}\}, \quad b_{ij} = \int_{S} \left(\lambda_1 \phi_j, \phi_i \right) dx 
+ \alpha \int_{\gamma^S} \left( \lambda_2 \hat{\phi}_j, \hat{\phi}_i \right) ds.
\end{split} 
\nonumber
\end{equation}
For definition of the multiscale basis functions, we select the first $M$ eigenvectors corresponding to the first $M$ smallest eigenvalues.

\subsection{Global formulation}
In the meshfree multiscale method the shape functions $W_i(x)$ defined in $S_i$ form the initial coarse space
\begin{equation} 
V^{\text{init}}_0 = \text{span} \{ W_i(x) \}_{i=1}^{N}.
\nonumber
\end{equation}
Here the shape functions $W_i(x)$ are defined as
\begin{equation}
W_i(x) = \frac{\phi_i(x)}{\sum_{j=1}^N \phi_j(x)},
\nonumber
\end{equation}
where $\phi_i(x)$ are kernel functions, which here are cubic splines
\begin{equation}
\phi(r) = 2 \left\{
\begin{array}{lr}
2/3 + 4 \ (r-1) \ r^2,& r \le 0.5,\\
4/3 \ (1-r)^3,& 0.5 \le r \le 1,\\
0,& 1 \le r,
\end{array}
\right.
\nonumber
\end{equation}
where $r$ is the normalized distance.

Accordingly, by multiplying the form function $W_i(x)$ by the eigenvectors $\psi_k^{\text{off}}$, we obtain the basis functions in the offline space $V_{\text{off}}$ 
\begin{equation}
\psi_{i, k} = W_i \psi_k^{ \text{off}}, \ 1 \le i \le N, \ 1 \le k \le M,
\nonumber
\end{equation}
where $k$ is an index of the eigenvector. Next, we construct the spectral multiscale space
\begin{equation}
V_{\text{off}} = \text{span} \{ \psi_{i, k}: 1 \le i \le N, \ 1 \le k \le M \}.
\nonumber
\end{equation}
Also, for further use, an operator matrix is constructed
\begin{equation}
R^T = [ \psi_1, \ldots, \psi_{N_c} ],
\nonumber
\end{equation}
where $N_c = N * M$ and $\psi_i$ are used to denote the nodal values of each basis function defined on the fine grid.

The transition matrix $R$ is then used to solve the following system on the coarse scale
\begin{equation}\label{eq:coarse_problem} 
M_\text{c} \frac{p_\text{c}^{n+1}-p_\text{c}^{n}}{\tau} + A_\text{c} p_\text{c}^{n+1} = 0,
\end{equation}
where 
\begin{equation}
M_\text{c} = R M R^T, \quad A_\text{c} = R A R^T.
\nonumber
\end{equation}

After solving the system \eqref{eq:coarse_problem}, we can go from a coarse-scale solution to a fine-scale solution using also the transition operator $R$ and the solution $p_\text{c}$
\begin{equation}
p_\text{ms} = R^T p_{c}, 
\nonumber
\end{equation}
where $p_\text{ms}$ is a fine-grid projection of the coarse-grid solution.

\begin{figure}[h!]
\begin{center}
\begin{minipage}[h]{0.32\linewidth}
\center{\includegraphics[width=\linewidth]{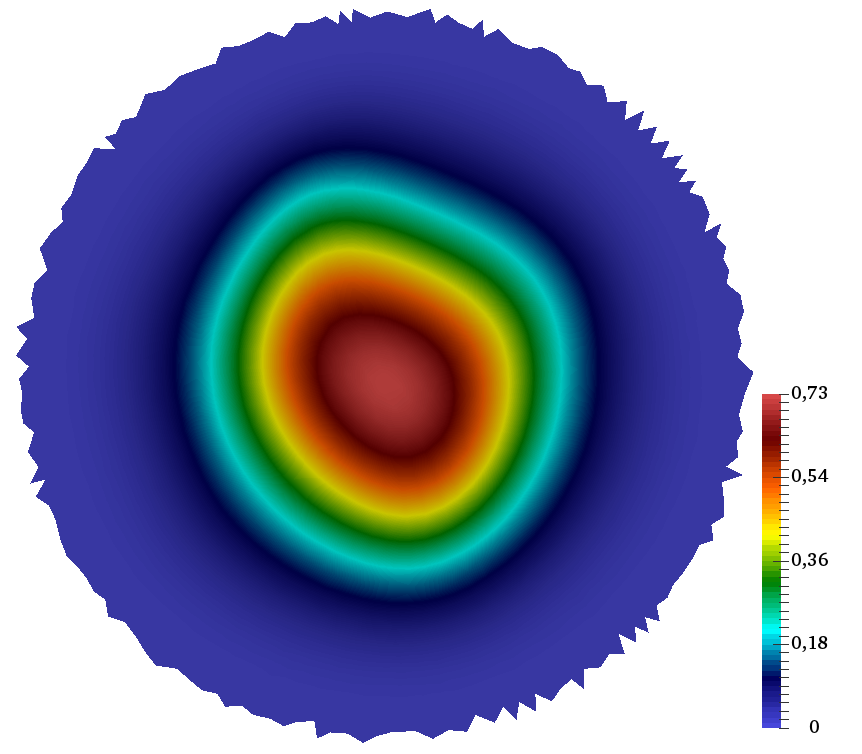}}
\end{minipage}
\begin{minipage}[h]{0.32\linewidth}
\center{\includegraphics[width=\linewidth]{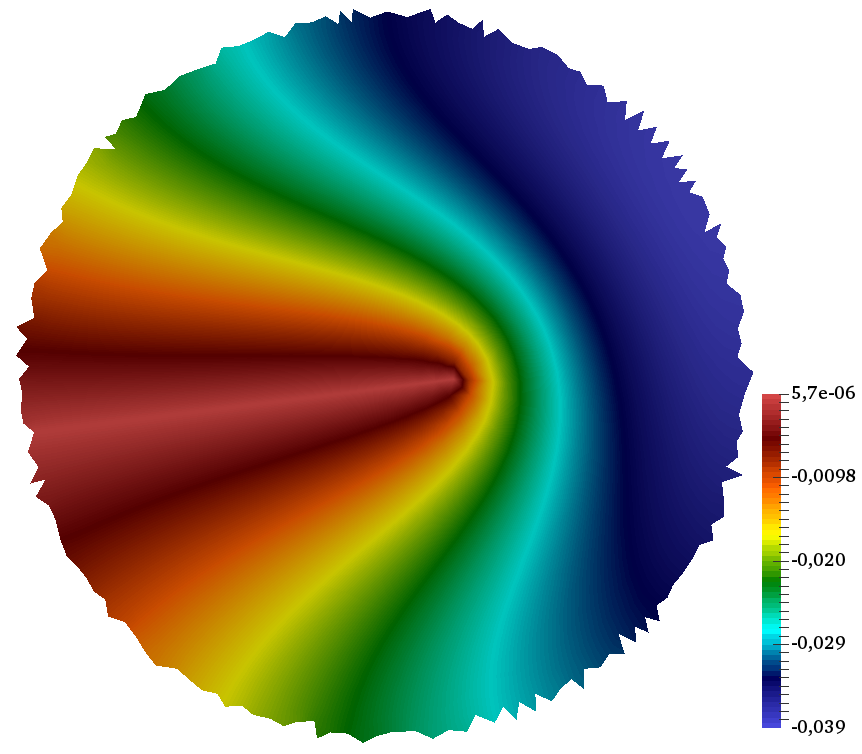}}
\end{minipage}
\begin{minipage}[h]{0.32\linewidth}
\center{\includegraphics[width=\linewidth]{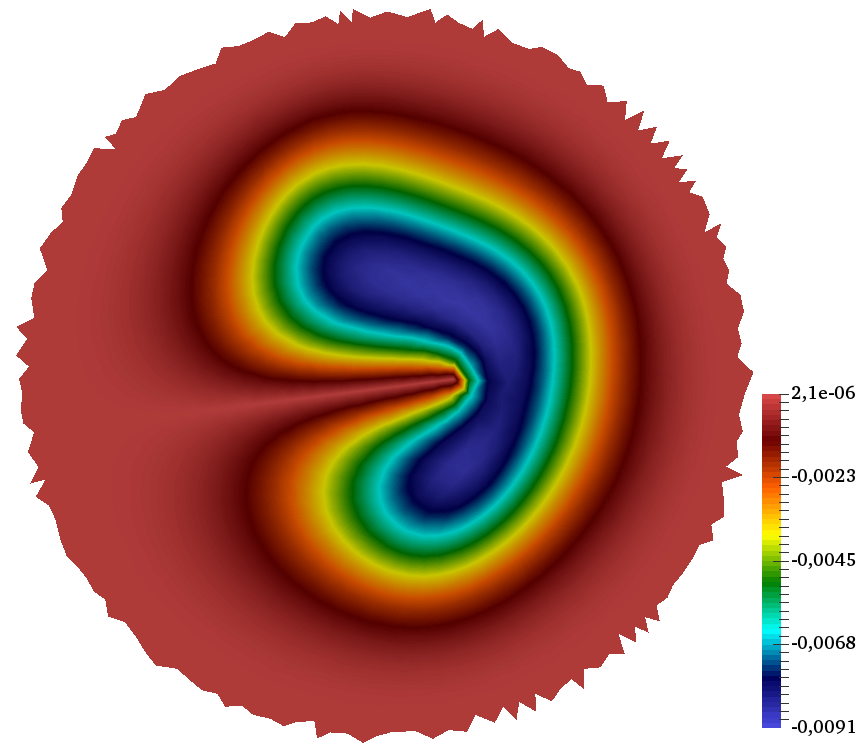}}
\end{minipage}
\end{center}
\caption{Illustration of multiscale functions. Left: shape function $W_i(x)$. Middle: eigenvector $\psi_k^{ \text{off}}$. Right: multiscale basis function $\psi_{i,k}$.}
\label{ris:basis1}
\end{figure}

\subsection{Partially explicit scheme}
Here we will discuss the splitting scheme at the coarse mesh nodes. All indices and functions correspond to the coarse mesh only. First we have the discrete multiscale function space $V_{\text{off}}$ and discrete solutions $p_c \in V_{\text{off}}$. We can represent the required function in the form of node values
\begin{equation}
p_c = \sum_{i=1}^N c_{i} \psi_{i},
\nonumber
\end{equation}
where for convenience, in contrast to the sum \eqref{eq:ms_sol}, we have reduced the index over the basis functions $k$, that is, $c_{i} = \sum_{j=1}^M c_i^j$ and $\psi_{i} = \sum_{j=1}^M \psi_i^j$. We use $p_c = (c_1, c_2, \ldots, c_N)^T$ to denote the vector of the required unknowns. 

\begin{figure}[h!]
\begin{center}
\center{\includegraphics[width=\linewidth]{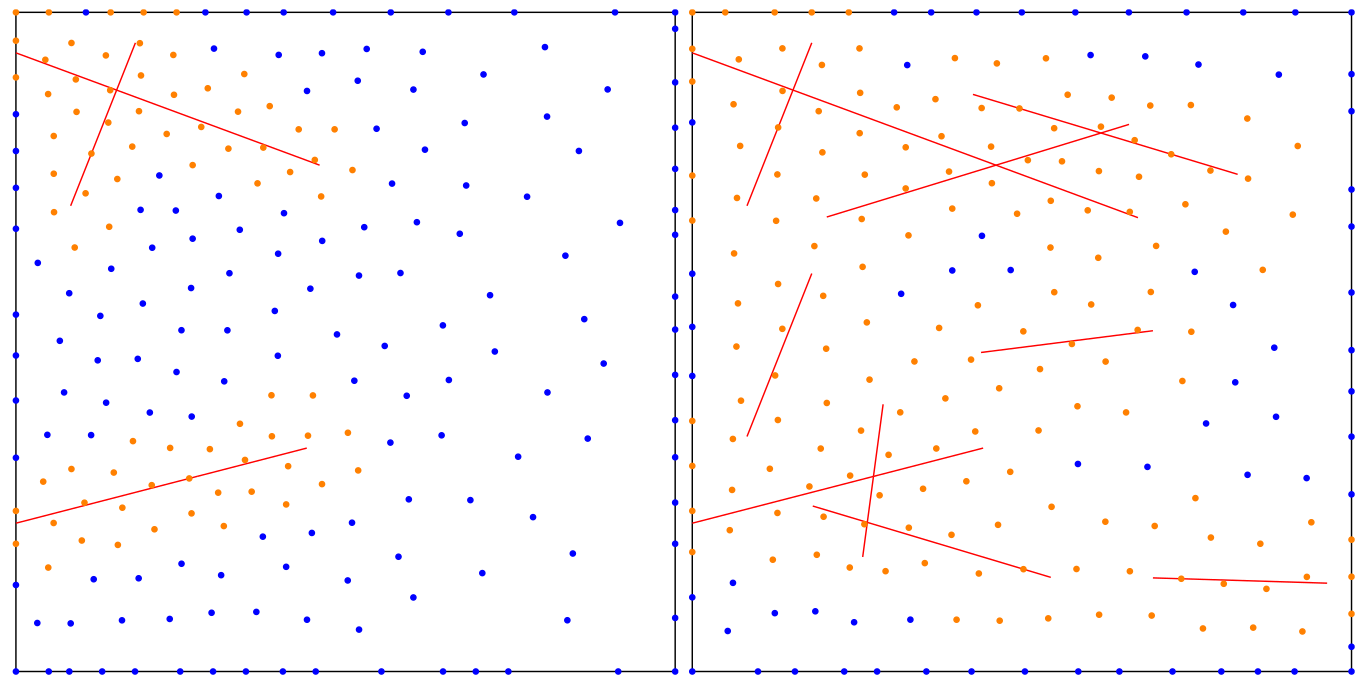}}
\end{center}
\caption{Point clouds with $N = 225$, where the orange dots are nodes for $p_I$ and the blue dots are nodes for $p_E$. Left: Test 1 with $N_I = 77$ and $N_E = 148$. Right: Test 2 with $N_I = 160$ and $N_E = 65$.}
\label{ris:pe_pic}
\end{figure}

Let $N = N_I + N_E$, where $N_I$ is the number of subdomains $S_i$ that contain fractures, and $N_E$ is the number without fractures. Then we can write
\begin{equation}
\sum_{i=1}^N c_i \psi_i = \sum_{i=1}^{N_I} c_i \psi_i + \sum_{i= N_I+1}^{N} c_i \psi_i,
\nonumber
\end{equation}
or $p_c = (p_I, p_E)^T$, where $p_I = (c_1, c_2, \ldots, c_{N_I})^T$ and $p_E = (c_{N_I+1}, c_{N_I+2}, \ldots, c_{N})^T$. Then, we have following system for the multiscale problem: find $p_c^{n+1} \in V_H$ such that
\begin{equation}
M_\text{c} \frac{p_\text{c}^{n+1}-p_\text{c}^{n}}{\tau} + A_\text{c} \begin{pmatrix}
p_I^{n+1} \\ p_E^{n} \end{pmatrix} = 0, \quad n = 0, \ldots , N_t - 1,
\nonumber
\end{equation}
where $p_I^{n+1}$ corresponds to implicit part of calculations, and $p_E^{n}$ to explicit ones. Thus we obtain a partially explicit scheme in time. The Figure \ref{ris:pe_pic} shows which nodes correspond to the solutions $p_I$ and $p_E$ for Test 1 and Test 2. 

\section{Numerical results}

This section presents numerical results for a single-phase filtration problem in a fractured porous medium, demonstrating that the proposed approach can predict an accurate approximation of the solution, independent of the contrast. The two-dimensional model problem defined in $\Omega = [0,80]\times[0,80]$ was considered. The fine meshes for Test 1 and Test 2 are similar and have 34,814 vertices and 68,983 elements and 34,949 vertices and 69,253 elements, respectively. For both tests in the multiscale method we use $N = 225$ and $M = 6$. For Test 1 we use $N_I = 77$ and $N_E = 148$ and for Test 2 we use $N_I = 160$ and $N_E = 65$.

The problem was solved numerically with the following values of the parameters: $c_m = 0.4$, \ $c_f = 1$, \ $k_m = 10^{-2}$, \ $k_f = 10^3$, \ $\alpha = 1$, \ $\mu = 1$, \ $p_0 = 1$, \ $g = 10$, \ $t_{max} = 900$, \ $N_t = 300$, \ $\tau = 3$. 

To approximate equations using the finite element method on a fine grid, construct systems and solve them, we use the open source computing platform FEniCS \cite{fenics} (LGPLv3).

\begin{figure}[h!]
\begin{center}
\begin{minipage}[h]{0.32\linewidth}
\center{\includegraphics[width=\linewidth]{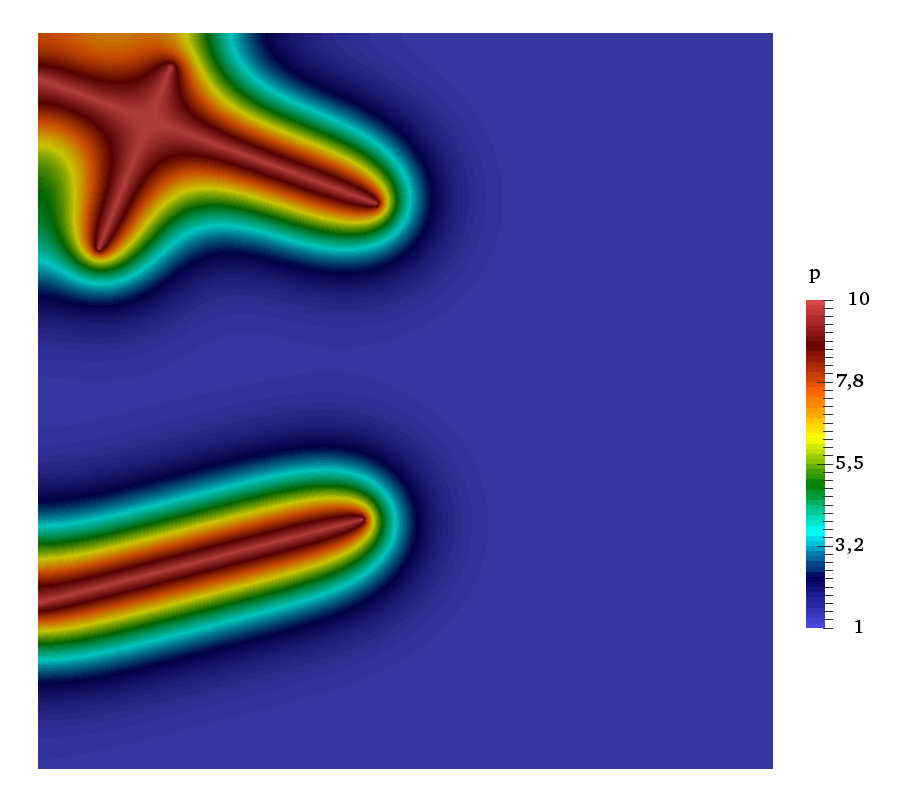}}
\end{minipage}
\begin{minipage}[h]{0.32\linewidth}
\center{\includegraphics[width=\linewidth]{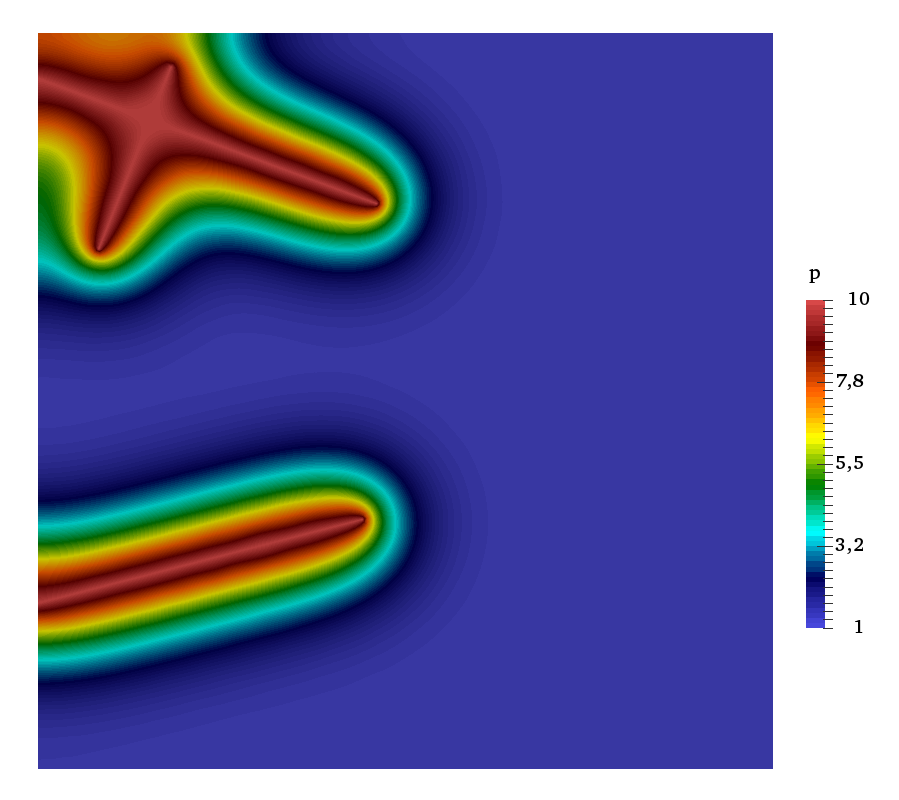}}
\end{minipage}
\begin{minipage}[h]{0.32\linewidth}
\center{\includegraphics[width=\linewidth]{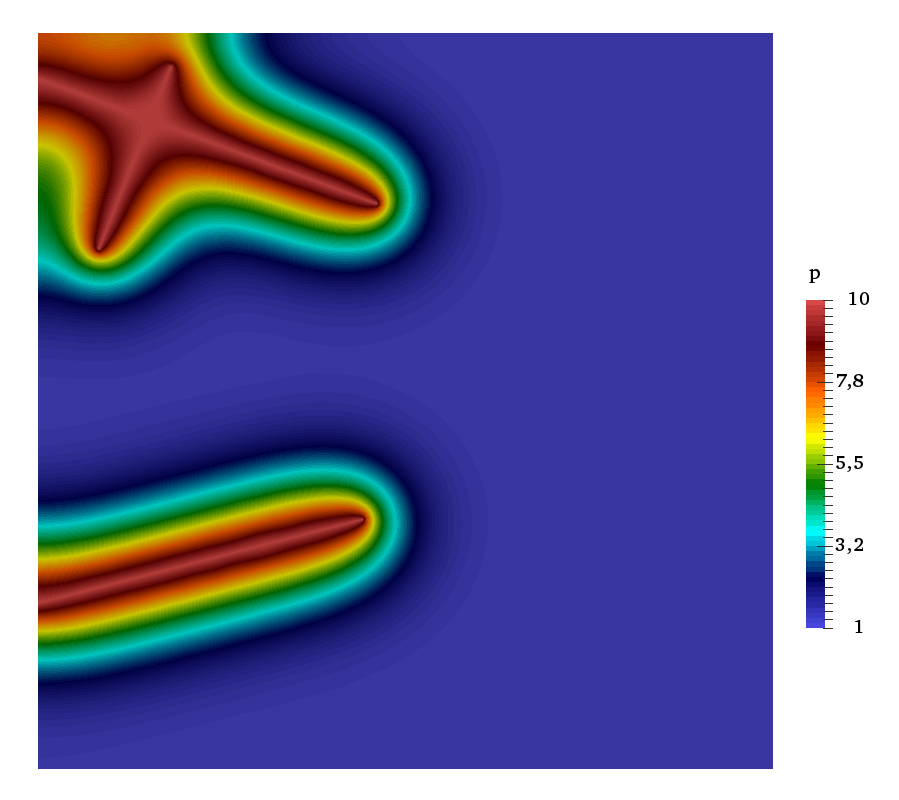}}
\end{minipage}
\end{center}
\caption{Solutions at the final time for Test 1. Left: reference solution. Middle: Meshfree GMsFEM with impicit method. Right: Meshfree GMsFEM with partially explicit method.}
\label{ris:sol1}
\end{figure}

\begin{figure}[h!]
\begin{center}
\begin{minipage}[h]{0.32\linewidth}
\center{\includegraphics[width=\linewidth]{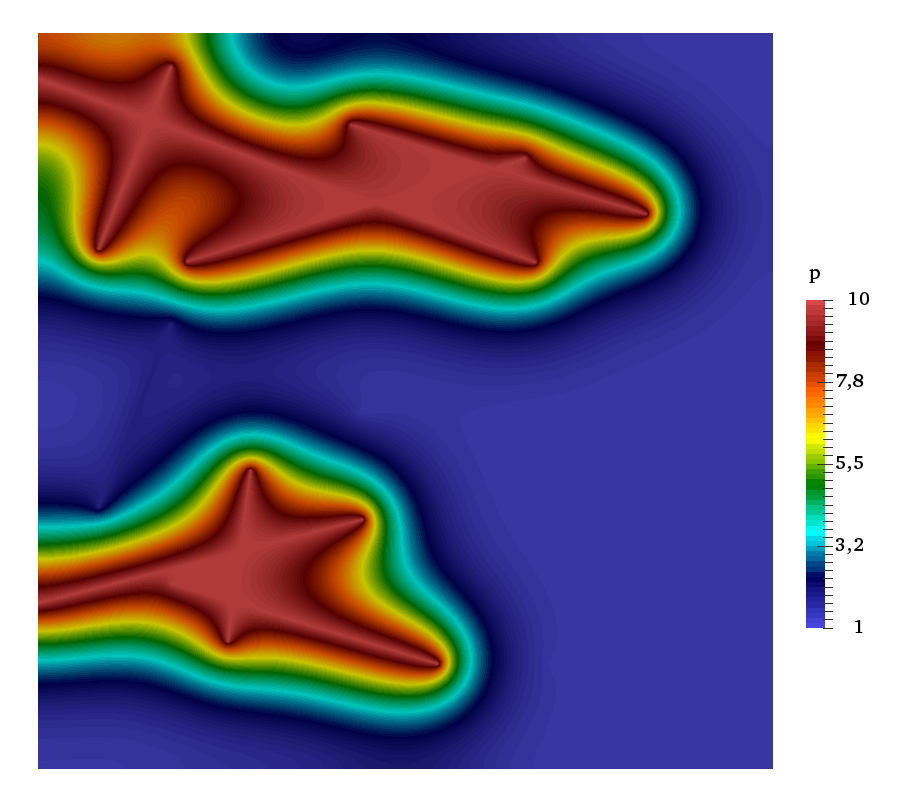}}
\end{minipage}
\begin{minipage}[h]{0.32\linewidth}
\center{\includegraphics[width=\linewidth]{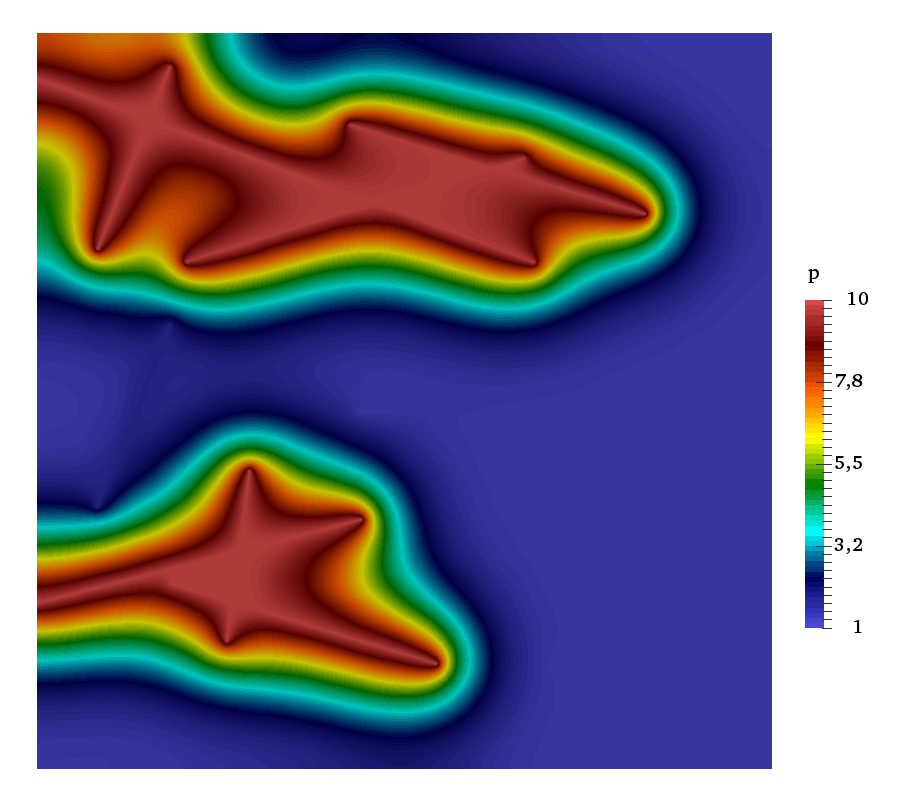}}
\end{minipage}
\begin{minipage}[h]{0.32\linewidth}
\center{\includegraphics[width=\linewidth]{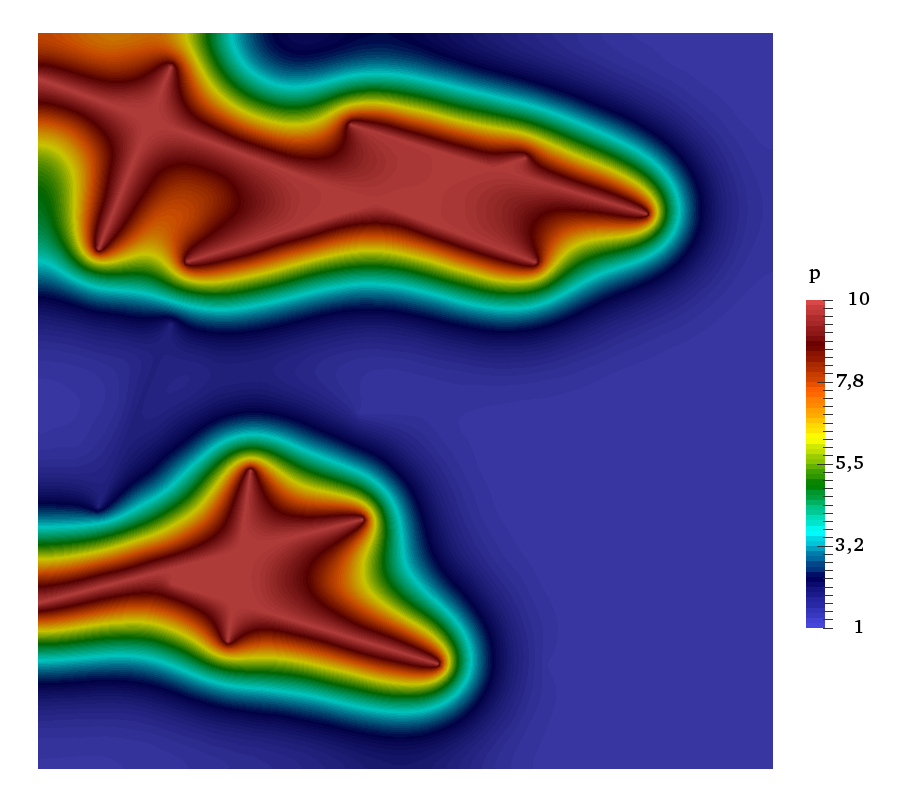}}
\end{minipage}
\end{center}
\caption{Solutions at the final time for Test 2. Left: reference solution. Middle: Meshfree GMsFEM with impicit method. Right: Meshfree GMsFEM with partially explicit method.}
\label{ris:sol2}
\end{figure}

Figures \ref{ris:sol1} and \ref{ris:sol2} show the solutions at the final time for Test 1 and Test 2, respectively. As we can see, the solutions obtained by multiscale methods with implicit and partially explicit discretization are similar to the reference solution on a fine grid.

To test our approach, we compare the solution obtained on a coarse grid with the solution obtained on a fine grid. The relative errors in $L^2$ and $H^1$ are calculated at each time as follows in percentages
\[ 
{||e||}_{L^2} = \frac{{||u_1-u_2||}_{L^2}}{{||u_1||}_{L^2}}\cdot100\%, \quad
{||e||}_{H^1} = \frac{{||u_1-u_2||}_{H^1}}{{||u_1||}_{H^1}}\cdot100\%,
\]
where ${||u||}_{L^2} = \int_{\Omega} u^2 dx, \ {||u||}_{H^1} = \int_{\Omega} (\grad u,\grad u) dx$, $u_1$ -- fine-grid solution and $u_2$ -- multiscale solution.

\begin{figure}[h!]
\begin{center}
\begin{minipage}[h]{0.49\linewidth}
\center{\includegraphics[width=\linewidth]{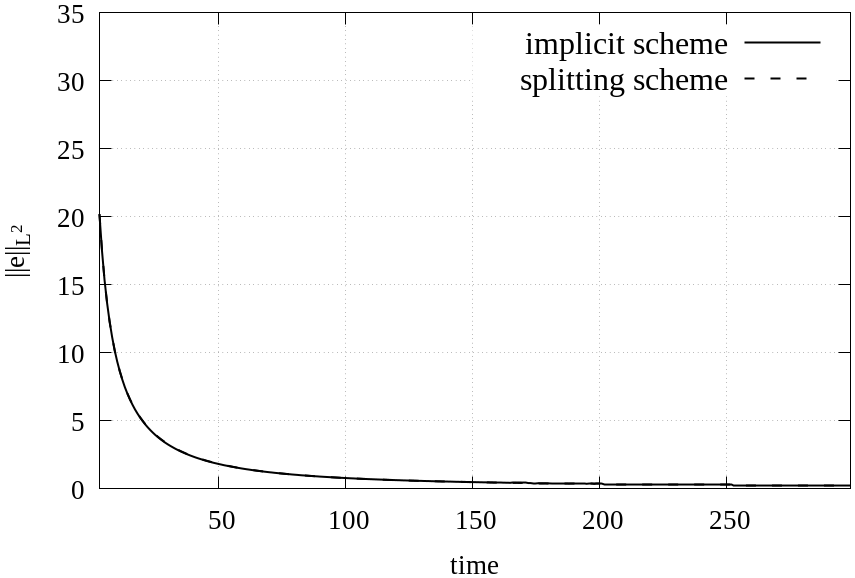}}
\end{minipage}
\begin{minipage}[h]{0.49\linewidth}
\center{\includegraphics[width=\linewidth]{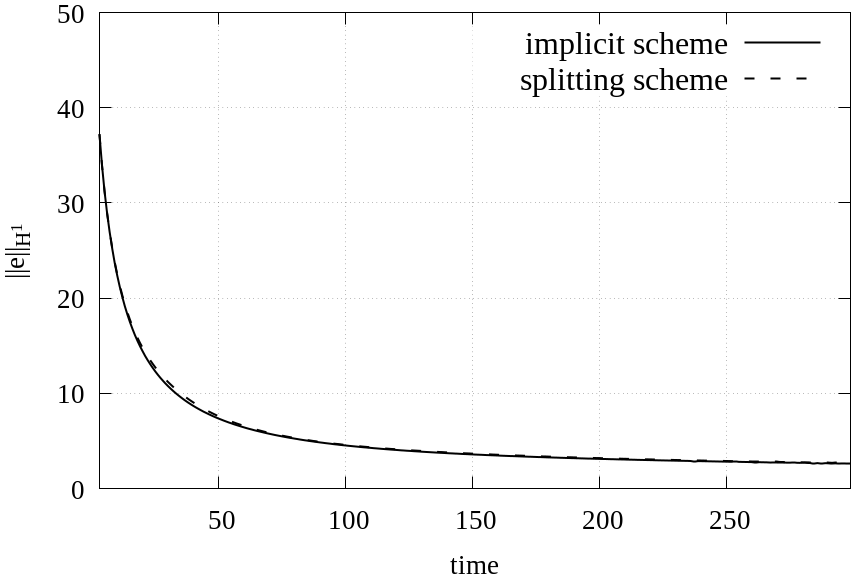}}
\end{minipage}
\end{center}
\caption{Plots of relative errors in $L^2$ (left) and $H^1$ (right) for Test 1 compared to the reference solution.}
\label{ris:error1}
\end{figure}

\begin{figure}[h!]
\begin{center}
\begin{minipage}[h]{0.49\linewidth}
\center{\includegraphics[width=\linewidth]{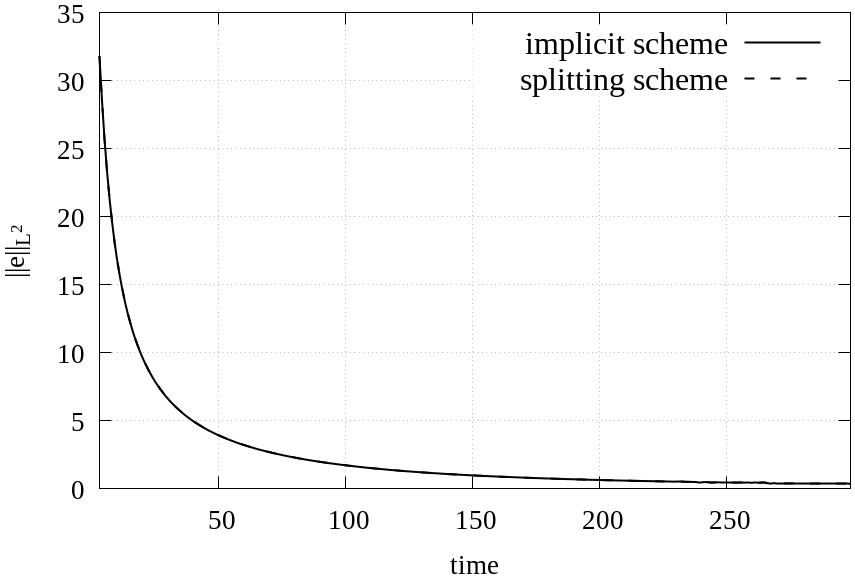}}
\end{minipage}
\begin{minipage}[h]{0.49\linewidth}
\center{\includegraphics[width=\linewidth]{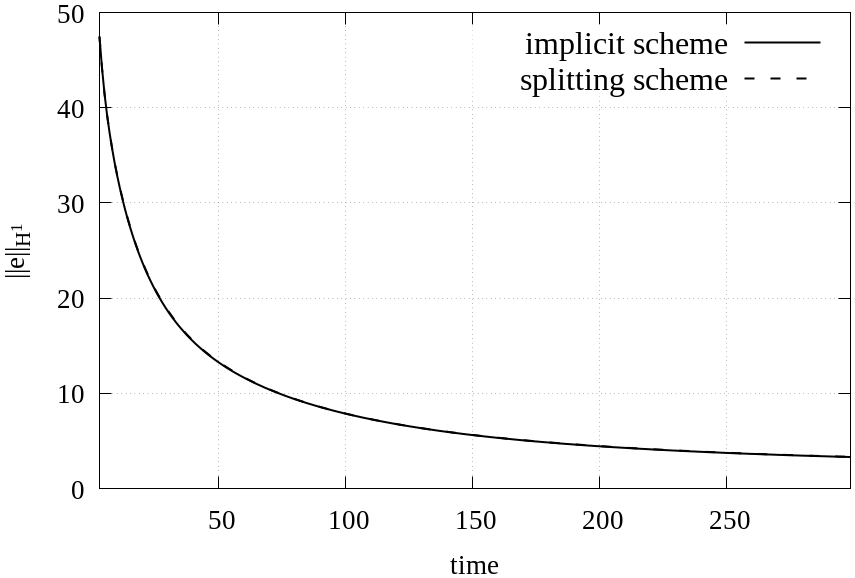}}
\end{minipage}
\end{center}
\caption{Plots of relative errors in $L^2$ (left) and $H^1$ (right) for Test 2 compared to the reference solution.}
\label{ris:error2}
\end{figure}

Figures \ref{ris:error1} and \ref{ris:error2} show the relative errors depending on the choice of time scheme for Test 1 and Test 2, respectively. Here, the purple curve represents the error caused by using the implicit time scheme, and the green curve represents the error of the partially explicit scheme. As we can see, the errors are almost the same. This indicates that our chosen partially explicit scheme provides reliable time discretization.

\begin{figure}[h!]
\begin{center}
\begin{minipage}[h]{0.49\linewidth}
\center{\includegraphics[width=\linewidth]{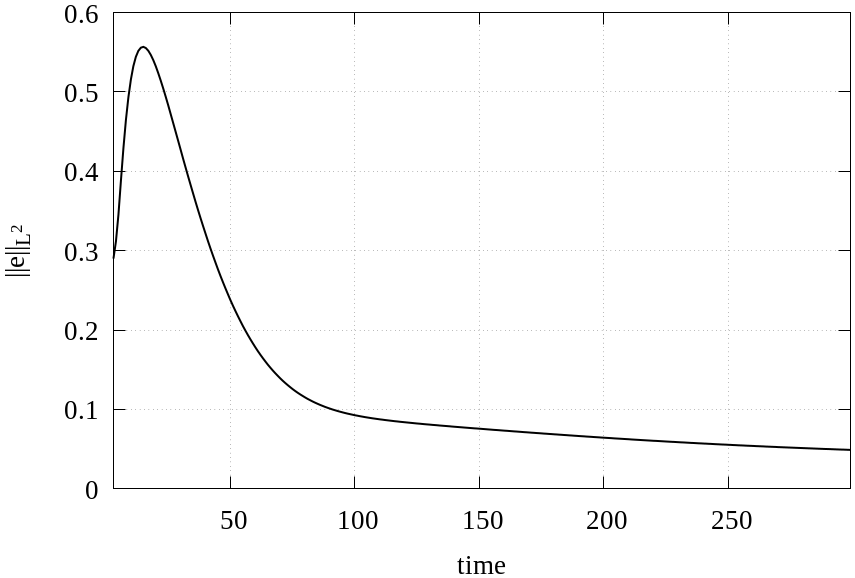}}
\end{minipage}
\begin{minipage}[h]{0.49\linewidth}
\center{\includegraphics[width=\linewidth]{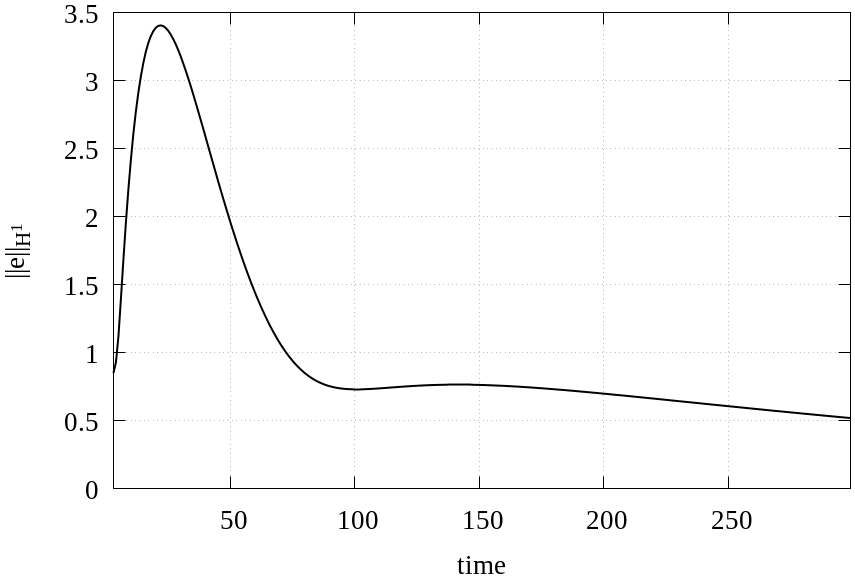}}
\end{minipage}
\end{center}
\caption{Plots of relative errors in $L^2$ (left) and $H^1$ (right) for Test 1 between impicit and partially explicit methods.}
\label{ris:error3}
\end{figure}

\begin{figure}[h!]
\begin{center}
\begin{minipage}[h]{0.49\linewidth}
\center{\includegraphics[width=\linewidth]{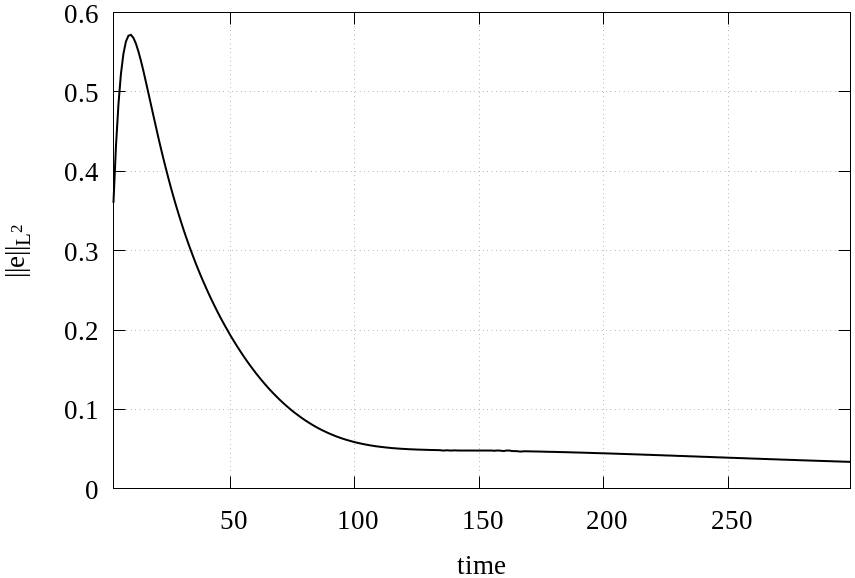}}
\end{minipage}
\begin{minipage}[h]{0.49\linewidth}
\center{\includegraphics[width=\linewidth]{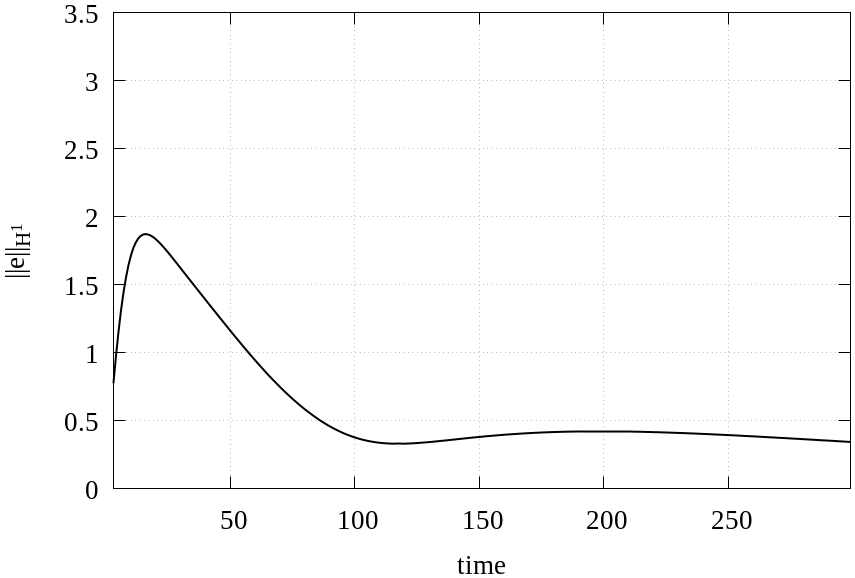}}
\end{minipage}
\end{center}
\caption{Plots of relative errors in $L^2$ (left) and $H^1$ (right) for Test 2 between implicit and partially explicit methods.}
\label{ris:error4}
\end{figure}

In Figures \ref{ris:error3} and \ref{ris:error4} we plot the relative errors of the solution obtained with the partially explicit scheme versus the solution with the implicit time scheme. Here we see that the errors in $L^2$ do not exceed $0.6\%$ for both tests, and in $H^1$ they do not exceed $3.5\%$ for Test 1 and $2\%$ for Test 2. Here, the number of nodes $N_I$ and $N_E$ for a partially explicit scheme have an influence, but it is insignificant, since for Test 1 we took twice as many nodes for the explicit part as for the implicit part and obtained good results.

\section{Conclusion}

In this paper, a time splitting scheme is proposed for the meshfree generalized multiscale finite element method. This scheme does not depend on the high contrast of the multiscale problem, but depends only on the permeability of the porous medium, which is much less than the permeability of fractures. It is shown numerically that our chosen partially explicit scheme provides reliable time discretization. The developed time discretization is beneficial for meshfree multiscale methods or other approaches that control the location of coarse mesh nodes to achieve high computational accuracy.


\begin{flushleft}
Djulustan Nikiforov,\\
North-Eastern Federal University, \\
NEFU, Belinsky St., 58, Yakutsk 677000, Russia,\\
Email: {\tt dju92@mail.ru},\\
\end{flushleft}


\begin{thebibliography}{99}
\bibitem{Efendiev1} {Y. Efendiev, T. Y. Hou}, {\it Multiscale finite element methods: theory and applications}, Springer Science \& Business Media, 4 (2009).
\bibitem{Efendiev2} {E. Chung, Y. Efendiev, T. Y. Hou}, {\it Multiscale Model Reduction: Multiscale Finite Element Methods and Their Generalizations}, Springer Nature, 212 (2023).
\bibitem{Nikiforov1} {S. Stepanov, D. Nikiforov, A. Grigorev}, {\it Multiscale Multiphysics Modeling of the Infiltration Process in the Permafrost}, Mathematics, 9.20 (2021), 2545.
\bibitem{Ammosov2} {D. Ammosov, M. Vasilyeva, E. T. Chung}, {\it Generalized Multiscale Finite Element Method for thermoporoelasticity problems in heterogeneous and fractured media}, Journal of Computational and Applied Mathematics, 407 (2022), 113995.
\bibitem{VPN1} {P. N. Vabishchevich}, {\it Additive Operator-Difference Schemes: Splitting Schemes}, Walter de Gruyter, (2013).
\bibitem{Ammosov1} {D. Ammosov, A. Grigorev, S. Stepanov, A. Tyrylgin}, {\it Partial learning using partially explicit discretization for multicontinuum/multiscale problems. Fractured poroelastic media simulation}, Journal of Computational and Applied Mathematics, 424 (2023), 115003.
\bibitem{Efendiev3} {E. T. Chung, Y. Efendiev, W. T. Leung, P. N. Vabishchevich}, {\it Contrast-independent partially explicit time discretizations for multiscale flow problems}, Journal of Computational Physics, 445 (2021), 110578.
\bibitem{Nikiforov3} {D. Nikiforov}, {\it Meshfree Generalized Multiscale Finite Element Method}, Journal of Computational Physics, 474 (2023), 111798.
\bibitem{Nikiforov4} {D. Y. Nikiforov, S. P. Stepanov}, {\it Modeling of Artificial Ground Freezing Using a Meshfree Multiscale Method}, Lobachevskii Journal of Mathematics, 44.3 (2023), 1206-1214.
\bibitem{dfm1} {M. V. Vasil'eva, V. I. Vasiliev, A. A. Krasnikov, D. Y. Nikiforov}, {\it Numerical simulation of single-phase fluid flow in fractured porous media}, Uchenye Zapiski Kazanskogo Universiteta. Seriya Fiziko-Matematicheskie Nauki, 159.1 (2017), 100-115.
\bibitem{dfm2} {V. I. Vasil'ev, et al.}, {\it Numerical solution of a fluid filtration problem in a fractured medium by using the domain decomposition method}, Journal of Applied and Industrial Mathematics, 12 (2018), 785-796.
\bibitem{cvt1} {L. Ju, Q. Du, M. Gunzburger}, {\it Probabilistic methods for centroidal Voronoi tessellations and their parallel implementations}, Parallel Computing, 28.10 (2002), 1477-1500.
\bibitem{cvt2} { M. Griebel, M. A. Schweitzer}, {\it A particle-partition of unity method for the solution of elliptic, parabolic, and hyperbolic PDEs}, SIAM Journal on Scientific Computing, 22.3 (2000), 853-890.
\bibitem{fenics} {A. Logg, K. A. Mardal, G. Wells (ed.)}, {\it Automated solution of differential equations by the finite element method: The FEniCS book}, Springer Science \& Business Media, 84 (2012).
\end{thebibliography}
\end{document}